\documentclass{amsart}
\usepackage{amssymb, amsfonts, amsmath,epsfig}
\usepackage{comment}
\begin{document}

\newtheorem{izr}{Theorem}[section]
\newtheorem{lm}[izr]{Lemma}
\newtheorem{pos}[izr]{Corollary}
\newtheorem{trd}[izr]{Proposition}
\newtheorem{conj}[izr]{Conjecture}
\newtheorem{df}[izr]{Definition}
\newtheorem{op}[izr]{Remark}

\newcommand{\Nn}{\mathbb{N}}
\newcommand{\Rr}{\mathbb{R}}
\newcommand{\Cc}{\mathbb{C}}
\newcommand{\Pp}{\mathbb{P}}
\newcommand{\I}{\mathbb{I}}
\newcommand{\Zz}{\mathbb{Z}}

\newcommand{\dk}[1]{\noindent \bf Proof #1 \rm }
\newcommand{\dok}{\noindent {\bf Proof.} \rm }
\newcommand{\ra}{\rightarrow \:}
\newcommand{\Ra}{\Rightarrow \:}
\newcommand{\lra}{\longrightarrow \:}
\newcommand{\Lra}{\Longrightarrow \:}
\newcommand{\Llr}{\Longleftrightarrow \:}
\newcommand{\ol}[1]{\overline{#1}}
\newcommand{\ul}[1]{\underline{#1}}
\newcommand{\ve}{\varepsilon}
\newcommand{\psh}{plurisubharmonic }
\newcommand{\nbhd}{neighbourhood }
\newcommand{\va}{\alpha}
\newcommand{\vb}{\beta}
\newcommand{\vph}{\varphi}
\newcommand{\tvp}{\tilde{\varphi}}
\newcommand{\PP}{{\bf P}}
\newcommand{\1}{\v{c}}
\newcommand{\2}{\v{C}}
\newcommand{\3}{\v{s}}
\newcommand{\4}{\v{S}}
\newcommand{\5}{\v{z}}
\newcommand{\6}{\v{Z}}
\newcommand{\dbar}{\overline{\partial}}
\newcommand{\Reg}{{\rm Reg }}
\newcommand{\Sing}{{\rm Sing }}
\newcommand{\Ker}{{\rm Ker }}
\newcommand{\pr}{{\rm pr }}
\newcommand{\Exp}{{\rm Exp }}

\newcommand{\m}[1]{\mathcal #1}

\newcommand{\Coker}{{\rm Coker }}
\newcommand{\supp}{{\rm supp }}
\newcommand{\D}{{\rm D }}
\newcommand{\Adj}{\rm{Adj}}
\newcommand{\Ricci}{\rm{Ricci}}

\def\med{\mathop{\rm med}\nolimits}
\def\min{\mathop{\rm min}\nolimits}
\def\emdim{\mathop{\rm emdim}\nolimits}
\def\dim{\mathop{\rm dim}\nolimits}
\def\codim{\mathop{\rm codim}\nolimits}
\def\rank{\mathop{\rm rank}\nolimits}
\def\corg{\mathop{\rm cork}\nolimits}
\def\Tr{\mathop{\rm Tr}\nolimits}
\def\Id{\mathop{\rm Id}\nolimits}
\def\Diag{\mathop{\rm Diag}\nolimits}

\newcommand{\noi}{\noindent}

\markboth{J. Prezelj}
{Positivity of metrics}

\title{Positivity of metrics on conic neighbourhoods of $1$-convex submanifolds}
\author{Jasna Prezelj}
\address{Jasna Prezelj,  Faculty of Mathematics, Natural Sciences and Information Technologies, University of Primorska, Glagolja\v ska 8, SI-6000 Koper, Slovenia and}
\address{Faculty of Mathematics and Physics, Department of Mathematics,
University of Ljubljana, Jadranska 19, SI-1000 Ljubljana, Slovenia and}
\address{Institute of  Mathematics,  Physics and Mechanics, Jadranska 19, SI-1000 Ljubljana,
Slovenia}
\email{jasna.prezelj@fmf.uni-lj.si}
\thanks{The author was supported by research program P1-0291 and by research project J1-5432 at Slovenian Research Agency. The part of this paper was written while the author was visiting the NTNU, Trondheim, Norway.}

\bigskip\bigskip\rm


\maketitle

\begin{abstract}
Let $\pi: Z \ra X$ be a holomorphic submersion from a complex manifold $Z$ onto a  $1$-convex manifold $X$ with  exceptional set $S$ and
  $a : X \ra Z$ a holomorphic section.  Let $\vph : X \ra [0,\infty)$ be  a plurisubharmonic exhaustion function which is strictly plurisubharmonic on $X \setminus S$ with  $\vph^{-1}(0) = S.$ For every holomorphic vector bundle $E \ra Z$ there exists a  neighbourhood $V$ of  $a(U\setminus S)$ for $U  = \vph^{-1}([0,c)),$ conic along $a(S),$ such that $E|_V$ can be endowed with Nakano strictly positive Hermitian metric.

   Let $g : X \ra \Cc, $ $g^{-1}(0) \supset S$ be a given holomorphic function. There exist finitely many bounded holomorphic vector fields defined on  a Stein neighbourhood ${V}$ of $a(\ol{U}\setminus g^{-1}(0)),$ conic along $a(g^{-1}(0))$ with zeroes of arbitrary high order on $a(g^{-1}(0))$ and such that they generate $\ker D\pi|_{a(\ol{U} \setminus g^{-1}(0))}.$ Moreover,  there exists a smaller neighbourhood ${V}' \subset V$ such that their flows remain in $\tilde{V}$ for sufficiently small times thus generating a local dominating spray.
\end{abstract}

\section{Introduction and  main theorems}

\noindent The main results of the present paper are  theorems \ref{main thm1} and \ref{main thm2}.

\begin{izr}[Nakano positive metric] \label{main thm1}
  Let $Z$ be an $n$-dimensional complex manifold, $X$  a $1$-convex manifold, $S \subset X$ its exceptional set,
  $\pi : Z \ra X$ a holomorphic submersion, $\sigma: E \ra Z$ a holomorphic vector bundle and $a: X \ra Z$
  a holomorphic section. Let $\vph : X \ra [0,\infty)$ be  a plurisubharmonic exhaustion function which is strictly plurisubharmonic on $X \setminus S$ and $\vph^{-1}(0) = S.$ Let $U  = \vph^{-1}([0,c))$ for some $c > 0$ be a given holomorphically convex set.
  There exist a neighbourhood $V_T$ of $a(\ol{U})$ in $Z$ and a Hermitian metric $h$ defined on  $E_{V_T \setminus \pi^{-1}(S)},$ such that
  \begin{trivlist}
     \item[] (a) $h$ has polynomial poles on $\pi^{-1}(S),$
     \item[] (b) there exists an open neighbourhood $V \subset V_T$ of $a(U\setminus S)$ conic along $a(S)$ such that $h$ is a
     Nakano positive Hermitian metric on  $E|_V,$
     \item[] (c) the curvature tensor $i\Theta(E)|_V$ has polynomial poles on $a(S)$ and is smooth up to the boundary elsewhere.
  \end{trivlist}
\end{izr}

\begin{izr}[Vertical sprays on conic neighbourhoods] \label{main thm2} With the same notation as above, let $g : X \ra \Cc$ be a holomorphic function with the zero set $N(g):= g^{-1}(0) \supset S$  and let $U =\vph^{-1}[0,c),$ $K \subset U, K \cap N(g) = \emptyset.$   There exist a Stein neighbourhood $V \subset Z$ of $a(U\setminus N(g))$ conic along $a(N(g))$ and finitely many bounded holomorphic vector fields $v_i$ generating $VT(Z) = \ker D\pi$ over a neighbourhood of $a(K)$ with zeroes
on $a(N(g))$ of arbitrarily high order. Consequently there exist  $\ve > 0$ such that the flows of $v_i$-s starting in a smaller conic neighbourhood $V' \subset V$ remain in $V$ for times $|t| <  \ve$ thus generating a local spray.
\end{izr}

The motivation for the present work was the paper \cite{pre} about the h-principle on $1$-convex spaces. Recall that a complex space $X$ is $1$-convex if it possesses a plurisubharmonic exhaustion function which is strictly \psh outside a compact set. There exists  a maximal nowhere discrete compact analytic subset $S$ of $X$ called the exceptional set.

In the proof we need a way of linearizing small perturbations of a given continuous section $a : X \ra Z,$ holomorphic on a given holomorphically convex open set $U,$ which are kept fixed on the exceptional set $S$ and are holomorphic on $U.$ This is usually done by using holomorphic sprays, i.e.  maps $s : U \times B^n(0,\ve) \ra Z,$ generated by holomorphic vector fields which span the vertical bundle $VT(Z) = \ker  D\pi$ on a neighbourhood $V \subset Z$ of $a(U)$ and are zero on $a(S).$
In the $1$-convex case such vector fields do not necessarily exist on the whole neighbourhood of $a(U)$ if $U$ intersects $S.$ In our application the condition on spanning $VT(Z)$ is needed on neighbourhood of  the set $a(K),$ $K \subset U$, where $K$ is a holomorphically convex compact set not intersecting $S;$ thus we can work with vector fields with zeroes (of high order) on $a(S)$ spanning $VT(Z)|_{a(K)}$ for $K$ satisfying $ K \cap S = \emptyset$ and it suffices if they are defined over a conic neighbourhood
of $a(U\setminus S).$ If they have zeroes of high enough order (with respect to the sharpness of the cone) their flows remain in the conic neighbourhood and thus generate the spray which dominates on a neighbourhood of $a(K).$ These vector fields are obtained as extensions of  vector fields defined on $a(X)$ which are  zero on a larger set, namely on the set $N(g) = g^{-1}(0),$ where $g : X \ra \Cc$ is a holomorphic function extended fibrewise constantly on $Z$ and such that $g(\pi^{-1}(S))=0$ and $N(g) \cap K = \emptyset.$
\begin{figure}[h!]
\begin{center}
 \epsfysize=45mm
 \epsfbox{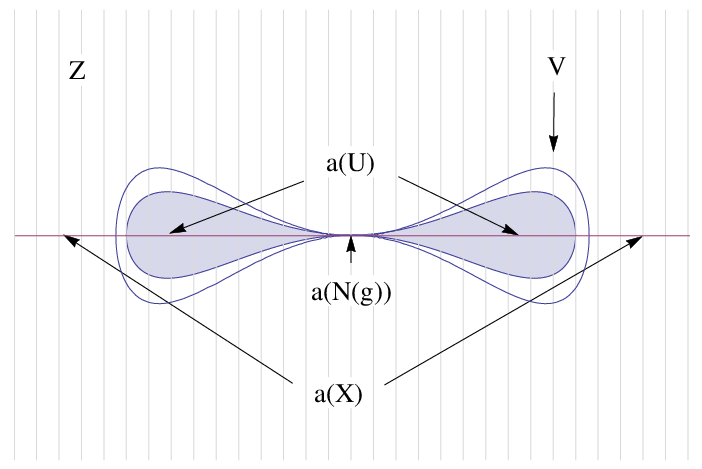}
 \caption{Conic neighbourhoods of $a(U \setminus N(g))$ in the submersion $Z \ra X$}
 \label{okolice}
\end{center}
  \end{figure}
Such extensions exist but it needs to be shown why they can be chosen to go to zero when approaching   $a(N(g)).$ This can be achieved by solving a suitable $\dbar$-equation with values in $VT(Z)$ and that is where we need the existence of Nakano positive metric. If $X$ were Stein the set $a(U)$ would have a basis of Stein neighbourhoods in $Z$ and a Nakano positive metric on $E|_V$ would be given by $h_E e^{-\psi}$ for some strictly \psh function $\psi.$ If $X$ is  $1$-convex then the set $a(U)$ does not necessarily have a basis of $1$-convex neighbourhoods   and on its neighbourhoods there are no  strictly \psh functions, since their Levi forms degenerate on exceptional sets. The construction of the metric and the sprays  is explained in the sequel.

\medskip

\noi{\bf Notation.}
The notation from the main theorems is fixed throughout the paper. Let $\omega_Z$ be a Hermitian $(1,1)$-form defined on the manifold $Z$ and $h_Z$ the corresponding Hermitian metric. Let $\sigma: E \ra Z$ be a holomorphic vector bundle of rank $r$  equipped with some Hermitian metric $h_E.$
The sets of the form $\pi^{-1}(U)$ are denoted by $Z_U.$ The local coordinate system in a neighbourhood $V_{z_0} \subset Z$ of  a point $z_0 \in a(U)$ is $(z,w),$ where $z$ denotes the horizontal and $w$ the vertical (or fibre) direction and $z_0 =(0,0).$ More precisely, every point in $a(U)$ has $w = 0$ and points in the same fibre have the same first coordinate.
The set $a(S)\cap V_{z_0} = \lbrace z_2 = 0, w = 0 \rbrace \cap V_{z_0},$ where $z_2$ is  a minimal set of generators of ${\m J}(S);$ if $S$ is a manifold, then they are just coordinate functions.
The function $\vph$ is extended to $Z$ fibrewise and we keep the same notation throughout the paper.
Its Levi form degenerates at most polynomially with respect to the distance from $Z_S.$ With the notation above this means that
the smallest eigenvalue of the Levi form does not go to zero faster than $\|z_2\|^{2k_0}$ for some $k_0 \in \Nn.$

The paper is organized as follows. In section $2$ we construct almost holomorphic global functions and the K\"{a}hler metric,  section $3$ is devoted to the proof of the first main theorem (theorem \ref{main thm1}), in section $4$ we solve $\dbar$-equation for $(n,q)$ and $(0,q)$-forms and in  section $5$ we prove the second main theorem (theorem \ref{main thm2}).

\section{Almost holomorphic global sections, \psh functions and the K\"{a}hler metric}

In this section we construct a K\"{a}hler metric on a conic neighbourhood of $a(U)$ using almost holomorphic global sections.
\begin{trd}[Almost holomorphic global sections]\label{sect}  Let $\sigma: E \ra Z$ be a holomorphic vector bundle. For
      every
      $l \in \Nn_0$ there exist a $k_l \in \Nn$ such that for  $ k \geq k_l$ there are finitely many  smooth sections $f_i$ of ${E},$ holomorphic in the vertical directions, such that they span ${ E}$ on some open neighbourhood $V_{T}$ of $a(\ol{U})$ in $Z$ except on $Z_{S}.$  Let $V_{z_0}\subset Z$ be a neighbourhood  of a point $z_0 \in a(U)$ such that $E|_{V_{z_0}}$ is trivial. Write $f_i = \sum f^{\lambda}_i e_{\lambda}$ with respect to some local frame $e_1,\ldots,e_r.$  If $z_0 \in a(S)$ there exists $C_i > 0$ such that for points $(z,w) \in V_{z_0}$ for sufficiently small $V_{z_0}$  we have
     \begin{eqnarray*}
         \|f_i(z,w)\| &\leq& C_1 \|z_2\|^k,  \\
         \|\dbar f_i(z,w)\| &\leq& C_2 \|w\|^{l+1}\|z_2\|^k,\\
         \|\partial f_i(z,w)\| &\leq& C_3 \|z_2\|^{k-1},\\
         \|\partial \dbar f_i\| &\leq& C_4 \|w\|^l\|z_2\|^{k-1}(\|w\| + \|z_2\|),\\
         \sum\|f^{\lambda}_i(z,w)\|^2 &\geq& C_5 \|z_2\|^{2k},\, \lambda = 1,\ldots,r.   \\
     \end{eqnarray*}
     If $z_0 \in a(U \setminus S)$ and $V_{z_0}$ is  sufficiently small we can replace $z_2$ by $1$ and obtain the estimates
     \begin{eqnarray*}
        \|\dbar f_i(z,w)\| &\leq& D_2 \|w\|^{l+1},\\
        \|\partial f_i(z,w)\| &\leq& D_3,\\
        \|\partial \dbar f_i\| &\leq& D_4 \|w\|^l,\\
        \sum\|f_i^{\lambda}(z,w)\|^2 &\geq& D_5, \lambda = 1,\ldots,r,\\
     \end{eqnarray*}
     for some $D_i > 0.$
\end{trd}

\op  Note that given $l$ the number $k$ can be chosen to be arbitrarily large.\\

\noi Before proceeding to the proof let us state a lemma on sections of  quotient sheaves.

\begin{lm}\label{dvig} Let and ${\mathcal E}$ be a coherent sheaf of sections of a holomorphic vector bundle $E \ra Z$ and denote by ${\mathcal Q} = {\mathcal J}(a(X))$ the ideal in ${\mathcal O}_Z$ generated by  the analytic set $a(X).$ Define
${\mathcal S}={\mathcal J}(a(S))^k ({\mathcal{E}}/{\mathcal Q}^{l+1})$ and let $F \in \Gamma(a(X),{\mathcal S})$
be a holomorphic section. Then for every point $z_0 \in a(S)$ there exist a local lift of $F_{z_0*}$ to a holomorphic section
$$
  F_{z_0}(z,w) = \sum_{|\va| = k, |\vb| \leq l, \lambda = 1,\ldots,r}  g_{\va\vb \lambda}(z){z_2^{\va}} w^{\vb}e_{\lambda} \in \Gamma(V_{z_0},E)
$$
in some local frame $\lbrace e_{\lambda} \rbrace$ and for $z_0 \in a(X \setminus S)$
there exist a local lift of the form $F_{z_0}(z,w) = \sum  g_{\va\vb \lambda}(z) w^{\vb}e_{\lambda} \in \Gamma(V_{z_0},E).$
\end{lm}

\dok The sheaf ${\mathcal S}$ is a finite dimensional vector bundle with coefficients in ${\mathcal J}(a(S))^k$ and it is supported on $a(X).$ Its sections represent Taylor series of vector fields in the $w$-variable up to order $l$ with coefficients in  ${\mathcal J}(a(S))^k.$ Since the statement is local we assume that $E$ is trivial and therefore it suffices to prove the result for functions.

 Let's assume that $z_0 \in a(S).$ In the given local coordinates near $z_0$ the generators of the ${\mathcal O}_Z/{\mathcal Q}^{l+1}$ are the germs $w^{\vb}_*$ ($\vb$ is a multiindex with $|\vb|\leq l$). Similarly, the minimal set of  generators of ${\mathcal J}(Z_S)^k$ is obtained from the minimal set of generators $z_2$ of the ideal ${\mathcal J}(Z_S)^k$, $\lbrace {z_2^{\va}}_*, |\va| = k, \va \in {\m A}\rbrace$, where we have omitted the unnecessary germs (if there are any) from the set of all monomials $\lbrace {z_2^{\va}}_*, |\va| = k \rbrace.$
  Their restrictions to $a(X)$ are the generators of ${\mathcal J}(a(S))^k.$
  Any element $G_{z_0*}$ of ${\mathcal S}_{z_0}$ is a finite sum of the form
  $G_{z_0*} = \sum {z_2^{\va}}_* (\sum g_{\va\vb *}  w^{\vb}_*),$ $g_{\va\vb*}\in {\mathcal O}_X.$
  Let $g_{\va\vb}$ be the local lifts of $g_{\va\vb *}$ to a neighbourhood of $z_0$ in $a(X)$ and fibrewise extended to $Z.$
  Then $G_{z_0}(z,w) = \sum  g_{\va\vb }(z) {z_2^{\va}} w^{\vb}$ is the desired lift defined on some neighbourhood  $V_{z_0}$ of $z_0.$
  For points $z \in a(X \setminus S)$ we replace the generators $z_2^{\va}$ by $1.$
\qed\\

\noi {\bf Proof of  proposition \ref{sect}.} By  theorem $A$ for relatively compact $1$-convex sets there exists $k_l \in \Nn$ such that for $k \geq k_l$
  there are finitely many sections $F_1,\ldots,F_m$  of  the sheaf ${\mathcal J}(a(S))^k ({\mathcal{E}}/{\mathcal Q}^{l+1})$  generating it on a neighbourhood  of $a(\ol{U})$ in $a(X).$

 Let $F$ be one of these sections and  $z_0 \in a(S).$ Choose a small product neighbourhood $V_{z_0}$ of $z_0$ in $Z$ with respect to the submersion $\pi: Z \ra X$ of the form $V_{z_0}=U_{z_0} \times B^{r_0}(0,\ve)$ in some local coordinates with $\pi \simeq \pr_1,$ the projection to the first coordinate.  By assumption $E$ is trivial on $V_{z_0}$ and the trivialization is given by the frame $e_1,\ldots,e_r.$
 Near $z_0$ the section $F$ has a local lift  $F_{z_0}$  defined on $V_{z_0}$ of the form
 $F_{z_0}(z,w) = \sum g_{\va\vb\lambda}(z) {z_2^{\va}} w^{\vb} e_{\lambda}$
  with coefficients as in lemma \ref{dvig}.
 Any other such lift for another choice of local generators $w$ coincides
 with this one up to order $l$  in  $w.$ If ${z_0}$ is not in $a(S)$ then we  assume that the closure of the neighbourhood $V_{z_0}$ does not intersect $Z_S.$
 Each $F_i$ thus defines an open covering of $a(\ol{U})$ in $Z$ and the latter has a locally finite subcovering.

  In the sequel we are examining the Taylor series of sections. They differ depending on the point $z_0$ under consideration.
  We focus on the case $z_0 \in a(S)$ and work in the usual coordinates $(z,w).$
  In the case $(z,0) \in a(U\setminus S)$  we replace the generators $z_2^{\va}$ of the ideal ${\mathcal J}(Z_S)^k$ in the estimates by the generator $1.$ Notice that there are constants $c_1, c_2 > 0$ such that $c_1 \|z_2\|^{2k} \leq \sum_{\va \in \m A} |z_2^{\va}|^2 \leq c_2 \|z_2\|^{2k}.$

  There exists a  locally finite product covering  $\lbrace V_j \cong U_j \times B^{r_0} \rbrace$ of $a(\ol{U})$ in $Z$
  by product neighbourhoods with respect to the submersion $Z \ra X$  finer than any of the above subcoverings. Let $\lbrace \chi_j \rbrace$
  be  a partition of unity subordinate to the product  covering which only depends on the base direction $z.$
  Summing up the local lifts $F_{ij}$ of $F_i$ on $V_j$ using this partition of unity
  we obtain sections  $ f_i(z,w)= \sum F_{ij}(z,w) \chi_j(z)$ on an open neighbourhood $U_Z$ of $a(\ol{U})$ in $Z$ which are holomorphic in the vertical direction and their nonholomorphicity is of the order  $\|z_2\|^k\|w\|^{l+1}$ as we see by expanding
  $F_{ij}$ in  Taylor series with respect to the vertical direction $w.$
  The terms in the expansion coincide up to order $l$ and therefore we have
  $F_{ij}(z,w) = F^l_i(z,w) + F_{ijl}(z,w),$ where $F_{ijl}$ are of order $\|z_2\|^{k}\|w\|^{l+1}$ and $F^l_i(z,w)$ have zeroes of order $k$ on $Z_S.$ Then $f_i(z,w) = F^l_i(z,w) + \sum F_{ijl}(z,w) \chi_j(z)$ and
\begin{eqnarray*}
  \|f_i(z,w)\| &\leq& C_1 \|z_2\|^k,\\
  \dbar f_i(z,w) &=& \sum F_{ij}(z,w) \dbar \chi_j(z) = \sum F_{ijl}(z,w) \dbar \chi(z)\\
  \partial f_i(z,w) &=& \sum \partial F_{ij}(z,w) \chi_j(z) + F_{ij}(z,w) \partial \chi_j(z),\\
  \partial \dbar f_i(z,w) &=& \sum \partial F_{ijl}(z,w) \wedge \dbar \chi_j(z) + F_{ijl}(z,w) \partial\dbar \chi_j(z).
\end{eqnarray*}
It is clear that there exist constants  $C_1-C_4$ and $D_2-D_4$ such that the claims hold.

 Because the sections generate $E$ on some neighbourhood of $a(\ol{U} \setminus S)$  the constant $D_5$ exists for a small neighbourhood of
 $z_0 \in a(\ol{U} \setminus {S})$ in $Z.$

  We still have to prove that the sections generate $E$ on a neighbourhood of $a(\ol{U})$ except on $Z_S$ to prove the existence of the constant $C_5.$ Since the statement is local, we may assume that $E$ is trivial, $E = V \times \Cc^r,$ with a local frame $e_1, \ldots, e_r.$
Let $A$ be the matrix with vector fields $f_i$-s as columns, $A = \lbrack f_1,\ldots,f_m\rbrack$ and consider the matrix $A A^*;$ they both have the same rank. We will show that the rank of $A$ equals $r$  by constructing a matrix $B = A G$  such that its columns will be approximately of the form $z_2^{\va} e_{\lambda}$  where $\va$ is a multiindex of order $k.$

By definition of $F_i$-s for any monomial $z_{2 *}^{\va}$ in ${\mathcal J}(a(S))^k$   at the point $z_0  \in a(S)$  there exist  coefficients $g_{\va i \lambda *}$ in the stalk ${\mathcal O}(a(X))_{z_0}$ such that
$F_{\va \lambda *} := \sum g_{\va i \lambda *} {F_i}_* = z_{2 *}^{\va }e_{\lambda}.$ Let $g_{\va i \lambda}$  be the functions on a neighbourhood $V_{z_0}$ of $z_0$ obtained by representing first the germs by  functions on a neighbourhood of $z_0$ in $a(X)$ and then  extending them fibrewise to  functions $g_{\va i \lambda}(z)$ depending only on $z.$
  Assume that the (local) sections $F_i$ of the sheaf are represented by sections of $E$ as above and denoted by the same letters. Then by definition of $F_i$-s we have
 $$
   F_{\va\lambda}(z,w) = \sum g_{\va i\lambda}(z) F_i(z,w) = z_2^{\va}e_{\lambda} + O(\|w\|^{l+1}\|z_2\|^k)
 $$
 and the same holds for the corresponding extensions $f_i,$ because they coincide with $F_i$-s to the order $l$ in $\|w\|,$
 $$
  F_{\va\lambda}(z,w) = \sum g_{\va i\lambda}(z) f_i(z,w) = z_2^{\va}e_{\lambda} + O(\|w\|^{l+1}\|z_2\|^k).
 $$
 Let $B$ be a matrix with $F_{\va\lambda}$ as columns for $\va \in {\m A},\lambda = 1,\ldots,r.$. We first write all $F_{\va\lambda}$ with $\lambda = 1$ and then with $\lambda = 2$ and so forth.
 Because the product $BB^*$ equals
 $$
   \sum |z_2^{ \va}|^2 I + O(\|w\|^{l+1}\|z_2\|^{2k}) = (\sum |z_2^{ \va}|^2)(I + O(\|w\|^{l+1}))
 $$
we conclude that
the vector fields $F_{\va \lambda}$ and therefore also the vector fields $f_i$ generate $E$ on a neighbourhood of $a(\ol{U})$ except on $Z_S.$ Since
$B = A G$ for the matrix $G$ defined by coefficients $g_{\va i \lambda}$ and because $\lbrace z_2^{\va}e_{\lambda} \rbrace $ is a subset of the canonical local generators $\lbrace z_2^{\va}w^{\vb}e_{\lambda} \rbrace$ the matrix
$G$ has full rank on a neighbourhood of $z_0.$ The matrix $B$ has full rank  on some open  neighbourhood $V_T$ of $a(\ol{U})$ except on $Z_S$ and so does $A.$ In other words, there exist a constant $C_5 > 0$ such that for every $\lambda$
$\sum \|f_i^{\lambda}(z,w)\|^2 \geq C_5\|z_2\|^{2k}$ provided $V_{z_0}$ is sufficiently small.
\qed\\

\begin{op}\label{lastne vrednosti}
Let $A \in \Cc^{r \times m}$ and $G\in \Cc^{m \times n}.$
Because $G$ has full rank at $z_0 = (0,0)$ in local coordinates, it has a singular value decomposition $G(z_0)= U^*\Sigma V$ with the matrix $\Sigma$ of  full rank equal to
$n.$ Then the $n \times n$ diagonal matrix $D = \Diag(d_1,\ldots,d_n)$ in $\Sigma$ is invertible. Since the singular values of
$AU^*$ and $BV^*$ are the same as those of $A$ and $B$ respectively we may assume that the matrices $U$ and $V$ are identities.
Then $G(z,w) = \Sigma + G_1(z,w), $ where $G_1(z,w) ={\m O}(\|(z,w)\|).$
Let $D(z,w) = D + {\m O}(\|(z,w)\|)$ be the upper $n \times n$ (invertible) square sub-matrix of $G.$
Then $A(z,w) G(z,w) D(z,w)^{-1} = A(z,w) \Sigma D^{-1} + \tilde{A}(z,w)$ with $\tilde{A}(z,w) = A{\m O}(\|(z,w)\|).$
Put $C(z,w):= A(z,w) \Sigma D^{-1} =A(z,w) I_{m,n} = B(z,w)D^{-1}(z,w) - \tilde A(z,w),$ where $I_{m,n}$ is the trivial inclusion $\Cc^n \hookrightarrow \Cc^m.$
Because of the properties of matrices $B$ and $D$ the matrix $CC^*$ is of the form
\begin{eqnarray*}
  CC^* &=& \Diag(\sum d_{i_1(\va)}^2|z_2^{\va}|^2,\ldots, \sum d_{i_r(\va)}^2|z_2^{\va}|^2)(I + {\mathcal O}(\|(z,w)\|)) + \\
  && \phantom{aaa} + {\mathcal O}(\|w\|^{l+1}\|z_2\|^{2k}) =\\
  &=& \Diag(\sum d_{i_1(\va)}^2|z_2^{\va}|^2,\ldots, \sum d_{i_r(\va)}^2|z_2^{\va}|^2)(I +{\mathcal O}(\|(z,w)\|)+ {\mathcal O}(\|w\|^{l+1}))
\end{eqnarray*}
so that its smallest eigenvalue decreases at most as $c_1 \sum|z_2^{\va}|^2$ and the largest is bounded from above
by $c_2 \sum|z_2^{\va}|^2.$
Then $A = \lbrack C| A_1\rbrack$ and since $AA^* = CC^* + A_1A_1^*$ the smallest eigenvalue of $AA^*$ does not decrease faster than $c_1\sum |z_2^{\va}|^2$ and because the entries of $A$ are bounded by $\|z_2\|^k$ the largest eigenvalue of $AA^*$
is bounded by $c_3\sum |z_2^{\va}|^2;$ the constants $c_1,c_2,c_3 $ are positive. Since all zeroes of the determinant  $\det (AA^*)|_{V_{z_0}}$ are on $V_{z_0} \cap Z_S$ it decreases polynomially with respect to $\|z_2\|$ on $V_T.$
\end{op}

\begin{op}\label{vext} Let $F$ be a holomorphic section of a holomorphic bundle $E$ over $Z$ defined on a neighbourhood of  $a(\ol{U})$ in $a(X).$  Then it is a section of
${\mathcal J}(a(S))^k ({\mathcal{E}}/{\mathcal Q}^{l+1})$ for $l = 0$ and some $k \geq 0.$ As in the proof of the above proposition there exist an almost holomorphic extension $f$ of  $F$ such that  $\dbar f(z,w)=  {\mathcal O}(\|z_2\|^k\|w\|).$
\end{op}

\subsection{Construction of a polynomially degenerating strictly plurisubharmonic function and the K\"{a}hler metric}
 In this section we describe the construction of a  function $\Phi$ which is strictly plurisubharmonic on a neighbourhood of $a(U\setminus S),$ conic along $a(S).$ Its Levi form decreases polynomially with the distance to $Z_S.$

With exactly the same construction as in the proposition \ref{sect} (we take a trivial line bundle) we produce a finite number of functions $\vph_{1,i}$ defined on an open neighbourhood $U'$ of $a(\ol{U})$ obtained from lifts of the sections of the sheaf \goodbreak
${\mathcal J}(a(S))^{k_1} ({\mathcal{J}}(a(U'))/{\mathcal J}^{l_1+1}(a(U')))$.
The sections are $0$ on $a(\ol{U})$, holomorphic to
order $l_1$ in the $w$-direction, have zeroes of order $k_1$ on $Z_S$ and such that away from $Z_S$ their vertical derivatives span the vertical cotangent bundle on a cone. The last assertion holds because near a point in $a(S)$ the functions are of the form
$$
  \vph_{1,i}(z,w) = \sum_{j,|\va|  = k_1}c_{ij\va}(z)z_2^{\va} w_j  + {\mathcal O}(\|z_2\|^{k_1}\|w\|^2)
$$
where $\lbrace z_2^{\va}, \va \in {\m A} \rbrace$ is a minimal set of generators of ${\mathcal J}(a(S))^{k_1}$ and $\va$ is  multiindex with $|\va| = k_1.$
Similarly as in the previous subsection we show that the functions $ z_2^{\va} w_j  $ for all possible $j,\va \in {\m A}$ are of the form
$z_2^{\va} w_j  = \sum g_{\va i j}(z) \vph_{1,i}(z,w) + {\mathcal O}(\|z_2\|^{k_1}\|w\|^2)$ and
$ z_2^{\va} dw_j = \sum g_{\va i j}(z) \partial_{w_j}\vph_{1,i}(z,w) + {\mathcal O}(\|z_2\|^{k_1}\|w\|).$
As before we conclude that  the forms $\partial_{w_j} \vph_{1,i}$ span the vertical cotangent bundle if $\|w\| \leq \|z_2\|$ and degenerate as $\|z_2\|^{k_1}.$  For points in $a(U \setminus S)$ with $\|z_2\|> \delta$  locally we have a uniform estimate, i.e. we replace $z_2$ by $1.$
Define $\vph_1 = \sum |\vph_{1,i}|^2$ whose Levi form
$$
  i\partial \dbar \vph_1 = i\sum \partial \vph_{i,1} \wedge \ol{\partial \vph_{i,1}} +  i\sum \ol{\dbar \vph_{i,1}} \wedge {\dbar \vph_{i,1}}
  +  i\sum \vph_{i,1} \ol{\dbar \partial \vph_{i,1}} + i\sum \partial \dbar  \vph_{i,1} \ol{\vph_{i,1}}
$$
has positive first two sums and all possibly negative terms are in the last two sums. Since they  involve at least one $\dbar \vph_{i,1}$ they
go to zero at least as $\|w\|^{l_1-1}.$  In the worst case the Levi form of
$$
  \Phi=\vph + \vph_1
$$
in coordinates $(z,w)$ is of the form
\begin{eqnarray*}
&& \begin{bmatrix}
     \|z_2\|^{2k_0} + \|w\|^2\|z_2\|^{2k_1-2} + \|w\|^{2l_1+2}\|z_2\|^{2k_1-2},& \|w\|\|z_2\|^{2k_1-1}\\
     \|w\|\|z_2\|^{2k_1-1},& \|z_2\|^{2k_1}
\end{bmatrix} \\
&&     +
\begin{bmatrix}
     \|w\|^{l_1+2}\|z_2\|^{2k_1-1},& \|w\|^{l_1}\|z_2\|^{2k_1}\\
     \|w\|^{l_1}\|z_2\|^{2k_1},& 0
\end{bmatrix},
\end{eqnarray*}
where the first matrix consists of the bound $\|z_2\|^{2k_0}$ for the smallest eigenvalue of the Levi form of $\vph$  and the first two terms of the above sum and is therefore positive and the second consists of the last two terms  and might be negative.
It is clear that  this form is positive on a neighbourhood of points from $a(U \setminus S).$ If we assume, say, that $\|w\| \leq \|z_2\|^{k_0 + 2}$ then the sum of such matrices is a positive definite matrix, since the diagonal block
$$
\begin{bmatrix}
     \|z_2\|^{2k_0} & 0\\
     0 & \|z_2\|^{2k_1}
\end{bmatrix}
$$
dominates. Instead of that we may  assume that $l_1 > 2 k_0$ and take the cone $\|w\| \leq \|z_2\|^2.$ In any case the Levi form $L\Phi$ is positive
on a conic neighbourhood of $a(U \setminus S)$ and the form
$$
  \omega = i\partial \dbar \Phi
$$
defines  the K\"{a}hler metric we are going to use.\\


\section{Nakano positive metric on a conic neighbourhood of a $1$-convex set}

In this section we first present some definitions and theorems on positivity of Hermitian metrics and then prove the first main theorem.

\subsection{Basic definitions and theorems on  positivity of Hermitian metrics}

We refer  to  Demailly's book Complex analytic and algebraic geometry \cite{dem} and
 recall some theorems from it.

Let $W$ be an $n$-dimensional complex manifold  and $E \ra W$ a holomorphic vector bundle equipped with a Hermitian metric $h.$ The matrix $H$ which corresponds to $h$ in a local frame $e_1,\ldots,e_r$ is given by
$$
  \langle u,v \rangle_h = \sum h_{\lambda \mu} u_{\lambda} \ol{v_{\mu}} = u^T H \ol{v}.
$$
The Chern curvature tensor $i\Theta(E)$ equals
$$
 i\Theta(E) = i\dbar(\ol{H}^{-1} \partial \ol{H}) = i\sum_{j,k}\Theta(E)_{jk} dz_j\wedge d \ol{z}_k.
$$
This can be considered as a matrix with $(1,1)$-forms as coefficients or as a $(1,1)$-form with matrices $i\Theta(E)_{jk}$ as coefficients.

If we denote the coefficient of $ dz_j\wedge d \ol{z}_k$  in the {\it column} $\lambda$ and the {\it row} $\mu$ by $c_{j k \lambda \mu},$ then
\begin{equation}\label{tensor}
 i\Theta(E) = i\sum c_{j k \lambda \mu} dz_j\wedge d \ol{z}_k \otimes e^*_{\lambda}\otimes e_{\mu},
\end{equation}
where  $1 \leq j,k \leq \dim W $ and $1 \leq \lambda, \mu \leq \rank E.$ Note that the matrix $\Theta(E)_{jk}$ has coefficients $\lbrace c_{jk\mu\lambda} \rbrace_{\mu,\lambda}.$
The bilinear form $\theta_E$ on $(TW \otimes E) \times (TW \otimes E)$  associated to $i\Theta(E)$ is defined by
\begin{eqnarray*}\label{forma0}
  \theta_E(u,v)  &=& \sum_{j,k} \langle \Theta(E)_{jk} u_j, v_k \rangle_h =
  \sum_{j,k} u_j^T \Theta(E)_{jk}^TH\ol{v}_k=
  \sum  c_{j k \lambda \mu} u_{j \lambda} \ol{v}_{k\nu} \langle e_{\mu},e_{\nu}\rangle_h = \\
  &=&\sum  c_{j k \lambda \mu} u_{j \lambda} \ol{v}_{k \nu} h_{\mu\nu},
\end{eqnarray*}
where $u = \sum_j ({\partial}/{\partial z_j}) \otimes u_j = \sum u_{j \lambda} ({\partial}/{\partial z_j}) \otimes e_{\lambda}$ and
$v= \sum_k ({\partial}/{\partial z_k})\otimes v_k= \sum v_{k \nu} ({\partial}/{\partial z_k}) \otimes e_{\nu}.$
In an orthonormal frame $e_1,\ldots e_r$ the form can be written as
\begin{equation}\label{forma}
  \theta_E = \sum  c_{j k \lambda \mu} (dz_j \otimes e^*_{\lambda}) \otimes \ol{(d{z}_k  \otimes e^*_{\mu})}.
\end{equation}

The form (\ref{forma}) gives rise to several positivity concepts. The `weakest'  is the Griffiths positivity which means that  the
form (\ref{forma}) is positive on the decomposable tensors $\tau = \xi \otimes v,$ $\xi \in TW,$ $v \in E$ so that
$$
  \theta_E(\tau,\tau) = \sum c_{j k \lambda \mu} \xi_{j} \ol{\xi}_{k} v_{\lambda} \ol{v}_{\mu}.
$$
The `strongest' is the Nakano positivity requiring that the form $\theta$  be positive on $\tau = \sum \tau_{j \lambda} ({\partial}/{\partial z_j}) \otimes e_{\lambda},$
$$
  \theta_E(\tau,\tau) = \sum c_{j k \lambda \mu} \tau_{j\lambda} \ol{\tau}_{k\mu}.
$$

In the case of holomorphic vector bundles the Griffiths curvature decreases in subbundles and increases in quotient bundles. This is not the case with Nakano positive bundles. Curvature in the sense of Nakano decreases in subbundles but does not increase in quotient bundles. In a related manner the dual of Nakano negative bundle is not necessarily Nakano positive. The connection between the two positivity concepts is described in the following

\begin{izr} [\rm Theorem VII-8.1, \cite{dem}] \label{dem3} If $E>_{\rm Grif} 0$ then $E\otimes(\det E) >_{\rm Nak}0.$
\end{izr}

Let $H$ be a matrix defining the metric $h$ on $E$ in a local frame $e_1,\ldots,e_r$ and let $H(z_0) = I.$ Then at $z_0$ the following hold:
\begin{eqnarray*}
  && \theta_{E \otimes (\det E)} = \theta_E + \Tr_E(\theta_E) \otimes h, \mbox{ where }\\
  && \Tr_E(\theta_E)(\xi,\xi)= \sum_{1 \leq \lambda \leq r} \theta_E(\xi \otimes e_{\lambda},\xi \otimes e_{\lambda}), \, \xi \in TW.
\end{eqnarray*}
This means that if $E$ is Griffiths positive then $\det E$ is positive. Let $e = e_1\wedge \ldots \wedge e_r$ and $\tau = \sum \tau_{j \lambda} ({\partial}/{\partial z_j}) \otimes e_{\lambda}.$ Then $\|e\|= 1$ and
\begin{equation}\label{forma1}
 \theta_{E \otimes (\det E)}(\tau \otimes e,\tau \otimes e) = (\sum c_{j k \lambda \mu}\tau_{j\lambda}\ol{\tau}_{k\mu}  +
 \sum c_{j k \lambda \lambda} \tau_{j\mu}\ol{\tau}_{k\mu})\|e\|^2.
\end{equation}
The last sum comes from the induced metric $\dbar \partial \log \det H$ on $\det E. $   In matrix form it is represented as
$(\dbar \partial \log \det H) \Id_E$ and the curvature of the tensor product is
$$
  i(\dbar (\ol{H}^{-1}\partial \ol{H}) + (\dbar \partial \log \det H) \Id_E)\otimes \Id_{\det E}.
$$


\subsection{Proof of theorem \ref{main thm1}.}

The Nakano positive Hermitian metric on $V$ is obtained from the induced metric on the quotient bundle of the trivial bundle. We first construct an almost Griffiths positive metric, correct it to a Griffiths positive one and then simulate the tensor product by the determinant bundle $\det E$ using a suitable weight to obtain an almost Nakano positive metric: we consider $E$ as  $E = (E \otimes \det E) \otimes(\det E)^*$ and choose a weight $\Phi_1$ in such a way that the line bundle $(\det E)^*$ with the metric $h_{\det E^*}e^{-\Phi_1}$ is almost Nakano positive  and in the last step correct this metric with another weight to make it Nakano positive.  In order to do this we have to have finitely many sections of $E$ spanning $E|_V$ which are holomorphic to a high degree. The form which defines the metric is defined on $V_T \setminus Z_S$ with polynomial poles
on $Z_S$ but fulfills the positivity requirements only on a conic set.

If we were given a Nakano positive metric on a neighbourhood of $a(S)$ then this construction would not be needed because the positivity  could be achieved by using a weight of the form $e^{-\Phi},$ where $\Phi$ is strictly plurisubharmonic on a neighbourhood of $a(\ol{U}\setminus S)$ conic along $a(S).$ In general we do not have such a metric.\\

\noi{\bf Proof.}
By proposition \ref{sect} there exist finitely many smooth
vector fields $f_1,\ldots, f_m$ on an open neighbourhood $V_{T}$ of $a(U),$ holomorphic to order $l$ in the vertical direction, and zero of order $k$ on $Z_S$ defining a surjective vector bundle homomorphism $f: U_Z \times \Cc^m \ra E|_{U_Z},$ where $U_Z = V_{T}\setminus Z_S.$   Thus the bundle $E|_{U_Z}$ can be given the metric of $\ker f^{\bot}.$ Consider the mapping $f$ in some local chart, denote by $r$ the rank of the bundle and let $(z,w)$ be the local coordinates as usual. Then the mapping  $f$ can be represented as a
$r \times m$ matrix $A$ with coefficients $f_{ij}$ which are holomorphic up to order $l$ in the vertical direction and therefore $\dbar A \approx\|w\|^l.$ The
linear mapping given by $A$ has the inverse $A^{-1} : E|_{U_Z} \ra \ker f^{\bot}.$ Then for $u,v \in E|_{U_Z}$
we have
$$
  \langle u,v \rangle_{h_0} := \langle A^{-1}u, A^{-1}v \rangle,
$$
where the right scalar product is the usual one on $\Cc^m.$ By definition the matrix $H_0 = \lbrace h_{0,ij} \rbrace$ associated with the $(1,1)$-form which defines the scalar product is
$$
  \langle u,v \rangle_{h_0} = \sum h_{0,ij} u_i \ol{v_j} = u^{\top} H_0 \ol{v} = u^{\top} {A^{-1}}^{\top} \ol{A^{-1}} \ol{v}
$$
and has poles on $Z_S.$
So
$$
    H_0 =  \ol{ {A^{-1}}^*{A^{-1}}}.
$$
The Nakano curvature tensor can be calculated by the formula
$$
  \Theta(E)_0 = \dbar (\ol{H_0}^{-1}\partial \ol{H_0}).
$$
Before continuing let us express $\ol{H_0}^{-1}$ by the matrix $A.$ Since away from $Z_S$ the matrix $A$ has full rank it has at every point $z_0 \in U_Z$ a singular value decomposition
$$
   A = V \Sigma U^*,
$$
where $V, U$ are unitary matrices and $\Sigma$ is a $r \times m$ matrix with all entries equal to $0$ except those on the diagonal, $d_1,\ldots, d_r,$ which are square roots of eigenvalues of $AA^*.$ The partial inverse $A^{-1}$ is then given by $U \Sigma^{-1} V^*,$ where $\Sigma^{-1}$ is $m \times r$ matrix with only diagonal elements $d_1^{-1}\geq \ldots \geq  d_r^{-1} > 0 $
nonzero. We have
$$
  { A^{-1}}^* A^{-1} = V {\Sigma^{-1}}^{\top} U^* U \Sigma^{-1} V^* =  V D^{-2} V^*,
$$
where $D$ is a diagonal matrix with diagonal $d_1,\ldots, d_r.$
By construction we have
$$
   A A^* = V \Sigma U^*U \Sigma^* V^* = V D^2 V^*
$$
and so
$$
  (A A^*)^{-1} = V D^{-2} V^* = {A^{-1}}^* A^{-1}.
$$
This means that
$$
   \ol{H_0} = (A A^*)^{-1}.
$$
For an invertible matrix $B$ we have $\partial B^{-1} = -B^{-1} \partial B B^{-1}.$ The curvature is
\begin{eqnarray*}
  \dbar (\ol{H_0}^{-1}\partial \ol{H_0}) &=& -\dbar ((A A^*) (AA^*)^{-1} \partial (AA^*) (AA^*)^{-1})\\
  &=& -\dbar (\partial (AA^*) (AA^*)^{-1})\\
  &=& - \dbar\partial (AA^*) (AA^*)^{-1} + \partial (AA^*)\wedge\dbar(AA^*)^{-1}\\
  &=& - \dbar\partial (AA^*) (AA^*)^{-1} - \partial (AA^*) (AA^*)^{-1}\wedge\dbar(AA^*)(AA^*)^{-1}.
\end{eqnarray*}

We are interested in calculating the curvature tensor at some point $z_0.$ Let's make a change of coordinates such that
 $D(z_0) = I.$ Then $AA^*(z_0) = I$ and the above expression simplifies to
$$
  \dbar (\ol{H_0}^{-1}\partial \ol{H_0}) =   - \dbar\partial (AA^*)  - \partial (AA^*)\wedge \dbar(AA^*).
$$
Let us calculate each of the terms separately. The first one is
\begin{eqnarray*}
  \dbar\partial (AA^*) &=& \dbar ((\partial A) A^* + A (\dbar A)^*) \\
  &=& (\dbar\partial A) A^* - \partial A \wedge(\partial A)^* + \dbar A\wedge (\dbar A)^*  + A (\partial \dbar A)^*,
\end{eqnarray*}
and the second one is
$$
  \partial (AA^*) \wedge \dbar(AA^*) = ((\partial A) A^* + A(\dbar A)^*)\wedge((\dbar A) A^*+ A (\partial A)^*).
$$
All of the terms  containing $\dbar A$ are small when close to the  section $a(U).$ If $z_0 \in a(U\setminus  S)$ then they are $0$.
We  divide the curvature form into two forms: the one without the $\dbar A$ expressions is denoted by
$\Theta_1$ and the remaining part by $\Theta_2.$
Then
$$
  \Theta_1 = -(- \partial A \wedge (\partial A)^*) -\partial A A^*\wedge A (\partial A)^* =\partial A \wedge (\partial A)^* -\partial A (A^* A)\wedge (\partial A)^*.
$$
Denote by $A_s$ the $s$-th column of $A.$ Since we have chosen $D(z_0) = I$ we have $A^*A = \pr_{\Cc^r}$ and this means that
$$
  \Theta_1(\xi \otimes v,\xi \otimes v) = \sum_{s = 1}^m  |\langle \partial A_s (\xi),v \rangle |^2 - \sum_{s = 1}^r  |\langle \partial A_s (\xi),v \rangle |^2  \geq 0
$$ is nonnegative on $V_T \setminus Z_S.$

If we multiply our initial trivial metric by $e^{-\Phi}$ the curvature tensor gets an additional term
$L\Phi,$ where $L\Phi$ denotes the Levi form of $\Phi$ and thus the form becomes strictly positive on
$a(U \setminus S)$ and consequently the bundle has positive Griffiths curvature at least on some open neighbourhood of $a(U\setminus S).$ We claim that it can be chosen to be conic.

 Wherever $\Phi$ is  strictly plurisubharmonic we are adding a strictly positive
 $(1,1)$-form. The bad news is that $\Phi$ is such only on a conic neighbourhood and its Levi form decreases polynomially as we approach $Z_S.$ But if we manage to show that the form $\Theta_2$ goes to $0$ even faster, then we can make Griffiths
 curvature positive on a conic neighbourhood.  In order to find the rate of decreasing we must work in
 ambient coordinates (and hence can not assume that $D(z_0) = I$ if  $ z_0 \in a(S)$).
 The form $\Theta_2$ is therefore equal to
 \begin{eqnarray*}
   \Theta_2 &=& (-\dbar\partial A A^* - \dbar A\wedge (\dbar A)^*  - A (\partial \dbar A)^*)(AA^*)^{-1} + \\
   & & - (\partial A A^* + A(\dbar A)^*)(AA^*)^{-1}\wedge(\dbar A A^*)(AA^*)^{-1} + \\
   & &  - A(\dbar A)^*(AA^*)^{-1}\wedge A (\partial A)^*(AA^*)^{-1}.
 \end{eqnarray*}

By construction the $\det (AA^*) = 0$ only on fibres above $S$ and goes to $0$ polynomially with respect to distance from
the $Z_S.$ If $z$ denotes the horizontal directions we have
$\det (AA^*) \geq c \|z_2\|^{n_2}$ for some constant $n_2$ (by remark  \ref{lastne vrednosti} the constant is in fact $n_2 = 2rk$). Because of noninvertibility of $AA^*$ the form
 $\Theta_2$  has poles and they are hidden in the determinant $\det (AA^*).$ Each term involving $(AA^*)^{-1}$ also involves a term of the form
$\dbar A \approx \|w\|^{l+1}\|z_2\|^k.$ So if $\|w\| \leq c \|z_2\|^{n_2 + n_3}$ for some $n_3 \in \Nn$ all the terms  will go to $0$
at least as $\|z_2\|^{n_3}$ inside this cone as we approach the set $a(S).$ If $n_3$ is large enough the possible negativity of $\Theta_2$ will be
compensated by the Levi form $L\Phi.$  Since we only have Griffiths nonnegative curvature it can be made strict by adding another factor $e^{-\Phi}.$ The new (now Griffiths positive) Hermitian metric on $E$ is denoted by
$$
  h_1 = h_0 e^{-2\Phi}.
$$

\begin{op}  Let $i\Theta_i = i\sum \Theta(E)_{jk}^i dz_j \wedge d\ol{z}_k. $ We may assume that at a given point after a unitary change of coordinates
we have $L\Phi = \sum \sigma_j dz_j \wedge d\ol{z}_j$ where $\sigma_j \geq c\|z_2\|^{2\max(k_0,k_1)}.$
 Let the  bilinear form $\theta$ be associated to $\Theta$ in the metric $h_0$ and let $\theta^1$ be associated to $\Theta^1  = \Theta + 2 L\Phi \Id_E$
 in the metric $h_1.$ The quadratic form for Griffiths curvature is
 \begin{eqnarray*}
   \theta^1(\xi \otimes v,\xi \otimes v) &=& \left(\sum \xi_j \ol{\xi}_k v^T \Theta(E)_{jk}^{1T} H_0  \ol{v} + \sum \xi_j \ol{\xi}_k v^T \Theta(E)_{jk}^{2T} H_0  \ol{v} + \right. \\
   && + 2 \left. \sum \sigma_j |\xi_j|^2 v^T H_0 \ol{v} \right) e^{-2\Phi}
 \end{eqnarray*}
for  $\xi \otimes v = \sum \xi_j v (\partial/\partial z_j).$
 The first form is nonnegative and the third degenerates in the worst case as $\|z_2\|^{2\max(k_0,k_1)- 2k}$  by remark \ref{lastne vrednosti}. The second form has coefficients
 bounded by  $\|z_2\|^{n_3 - 2k}$ when approaching $Z_S$ and for large $n_3$ they are smaller than  $\|z_2\|^{2\max(k_0,k_1)- 2k}$ and for an even larger $n_3$ they go to zero as $\|z_2\| \ra 0.$
\end{op}

Let $H_1$ be the matrix representing $h_1$ in a local frame of $E.$ Then the determinant bundle has
a metric given by $\tau_1 = \det (h_{1, \lambda \mu})$ and  since the curvature of $\det E$ is positive, we have
$$
  -\partial \dbar \log \tau_1 = \partial \dbar \log \tau_1^{-1}> 0.
$$
Consider the induced metric on the dual bundle $E^*.$ Let $e_1,\ldots, e_r$ be a local  frame  of $E$ and $e_1^*,\ldots e^*_r$ the dual frame. Each $e_{\lambda}^*$ can be
represented as the scalar product by the vector $f_{\lambda}$ satisfying the equation $\langle e_{\mu}, f_{\lambda}\rangle_{h_1} = \delta_{\lambda\mu}$ or $H_1 \ol{F}  = I$ where
 $F = \lbrack f_1,\ldots,f_r \rbrack .$ Then the induced scalar product is given by the matrix
$\ol{F^T H_1 \ol{F}}  = F^* = {H_1^T}^{-1}.$
  The induced metric $\det (h_1)^*$ on $\det E^*$ in the dual coordinates is thus represented by $\tau_1^{-1}.$  Let $v_1^*,\ldots,v_k^*$ be almost holomorphic sections of  $(\det E)^*$ given by  proposition \ref{sect}. They generate the bundle on a neighbourhood of $a(U)$ in $Z$ except over the fibres over $a(S).$ Then we can multiply the metric $h_1$ by the weight
$e^{-\log \Phi_1} = \Phi_1^{-1},$
$$
  \Phi_1 = \sum_i \langle v_i^*,v_i^*\rangle_{\det (h_1)^*},
$$
to obtain the metric
$$
  h_2 = h_1 e^{-\log \Phi_1}.
$$
In the local frame $e_1,\ldots e_r$ of $E$ we have with $e^* := (e_1\wedge\ldots \wedge e_r)^*$ the norm
$$
   \langle e^*,e^*\rangle_{ \det (h_1)^*} = \tau_1^{-1}
$$
and since $v_i^* = \va_i e^*$ for some almost holomorphic functions $\va_i$ we have
$$
   \langle v_i^*,v_i^*\rangle_{ \det (h_1)^*} = \tau_1^{-1}|\va_i|^2
$$
and so the weight equals
$$
  \Phi_1 = \sum (\tau_1^{-1}|\va_i|^2) = \tau_1^{-1}\sum |\va_i|^2.
$$
The metric is
$$
  h_2 = h_1 \tau_1 \frac{1}{\sum |\va_i|^2}
$$
and has again polynomial poles only on $Z_S$ if restricted to some small  neighbourhood of $a(\ol{U})$ in $Z.$
The curvature tensor is then
$$
  i \dbar(\ol{H_1}^{-1}\partial\ol{H_1}) + (i\partial \dbar \log \tau_1^{-1} + i\partial\dbar \log \sum |\va_i|^2)\Id_E
$$
and has polynomial poles on $Z_S.$

The first two terms represent the curvature tensor of $E \otimes(\det E) $;  it is Nakano positive by theorem (\ref{dem3}) wherever
$E$ is Griffiths positive. The last term would be  nonnegative if $\va_i$ were holomorphic. Since they are only almost holomorphic there may be negative terms hidden in the last sum of the curvature tensor. But all the negative terms are multiplied by terms of the form $\dbar \va_i$ and  only add terms that are bounded (and go to zero) on some conic neighbourhood:
\begin{eqnarray*}
  i\partial \dbar \log \sum \va_i \ol{\va_i} &=& \frac{i}{(\sum|\va_i|^2)^2}(\sum |\va_j|^2 \sum \partial \va_i \wedge \ol{\partial \va_i} - \sum \partial \va_j \ol{\va_j} \wedge \sum \va_i \ol{\partial \va_i}) + \\
  &&  - \frac{i}{(\sum|\va_i|^2)^2} \left(\sum \va_j \ol{\dbar \va_j}\wedge(\va_i \ol{\partial \va_i} + \ol{\va_i}\dbar{\va_i}) + \ol{\va_i \va_j} \partial \va_j \wedge \dbar \va_i \right) +\\
  && + \frac{i}{\sum|\va_i|^2}\left( \sum \va_i L \ol{ \va_i} + L\va_i \ol{\va_i} + \ol{\dbar \va_i} \wedge {\dbar \va_i}    \right).
\end{eqnarray*}
The first line is positive by the Lagrange identity and the rest is potentially negative. Take a point $(z,0) \in a(U\setminus S).$ There we have $\dbar \va_i = 0$ and $\dbar \va_i(z,w) \approx \|w\|^{l_2}$ for some  $ l_2 > 2$ otherwise. On a neighbourhood of $(z,0) \in a(S)$ we have for some $k_2> 2$ by proposition \ref{sect} the estimates
\begin{eqnarray*}
  &&\sum|\va_i|^2 \approx \|z_2\|^{2 k_2},\\
  &&\dbar \va_i(z,w) \approx \|w\|^{l_2}\|z_2\|^{k_2},\\
  &&\partial \va(z,w) \approx \|z_2\|^{k_2 - 1}\\
  &&L \va_i(z,w) \approx \|w\|^{l_2-1}\|z_2\|^{k_2-1}(\|z_2\| + \|w\|).
\end{eqnarray*}
So second and the third line of the Levi form  are  of the form
$$
  C_1\frac{\|w\|^{l_2}}{\|z_2\|} + C_2\|w\|^{2l_2} + C_3 \|w\|^{l_2} + C_4{\|w\|^{l_2 - 1}}{}
$$
and decrease polynomially in conic neighbourhoods of the form $\|w\| < \|z_2\|^{k_3}$ for $k_3$ large enough. Therefore in some conic neighbourhood
thin enough with respect to $\|w\|$ and sharp enough along $a(S)$ the negativity of these two terms can be compensated by the weight $e^{-C \Phi}$ for some positive constant $C$ as before. Since $S$ is compact there exist $C$ large enough for all $(z,0) \in a(S).$
The desired metric is therefore
$$
   h = h_3 = h_2 e^{-(C+1) \Phi} = h_0 e^{-((C+3) \Phi + \log \Phi_1)}, C > 0
$$
and has polynomial poles  on $Z_S,$ $h(z,w) \approx \|z_2\|^{-\kappa_1}h_E(z,w),$ $\kappa_1 \in \Nn.$
$\;$ \qed

\begin{op}\label{poli} Note that choosing a large $k_2$  produces a large pole on $Z_S$ in the weight. The  form $\theta_3$ corresponding to $h_3$ also has polynomial poles only on $Z_S.$
\end{op}


\section{$\dbar$-equation on conic neighbourhoods}

In this section we first present some results on $L^2$-methods on $\dbar$-equation from \cite{dem} and then solve the $\dbar$-equation for $(n,q)$ and $(p,0)$-forms.


\subsection{Basic theorems on $\dbar$-equation with values in a vector bundle}

Let $(W,\omega)$ be an $n$-dimensional K\"{a}hler manifold with the K\"{a}hler form $\omega = i \sum \gamma_i dz_i \wedge d\ol{z}_i$ , $E \ra W$ a vector bundle equipped with a Hermitian metric $h$ and let $H$ be the corresponding matrix in a local frame $e_1,\ldots,e_r.$
Let $i\Theta(E)$ be the Chern curvature tensor  and $\Lambda$ the adjoint of the operator
$u \ra u \wedge \omega$ defined on $(p,q)$-forms. The scalar product on $\Lambda^{p,q}(W,E)$ is   defined pointwise as
$$
  \langle  u_{JK\lambda}dz_J\wedge d\ol{z}_K \otimes e_{\lambda},  v_{J_1K_1\mu}dz_{J_1}\wedge d\ol{z}_{K_1} \otimes e_{\mu}\rangle =
  u_{JK\lambda}\ol{v_{JK\mu}}\gamma^{-J}\gamma^{-K} h_{\lambda \mu},
$$
if $J = J_1, K = K_1$ and $0$ otherwise;
$\gamma = (\gamma_1,\ldots,\gamma_n)$ and $J,K$ are multiindices, $|J| =|J_1| = p, |K|=|K_1|= q$.
Denote by $L^2_{p,q} (W,E)$ the space of $(p,q)$-forms with values in $E$ and with bounded $L^2$-norms with respect to the given metric $h$ and the form $\omega.$
Define the Hermitian operator $A_{E,\omega}$ as the commutator
$$
   A_{E,\omega} = [i\Theta(E), \Lambda].
$$

\begin{izr}[\rm Theorem VIII-4.5, \cite{dem}]\label{dem1} If $(W, \omega)$ is complete and $A_{E,\omega}> 0$ in bidegree $(p, q)$, then for any $\dbar$-closed form  $u\in L^2_{p,q} (W,E)$ with
$$
  \int_W \langle A^{-1}_{E,\omega} u, u \rangle dV < \infty
$$
there exists $v \in L^2_{p,q-1}(W,E)$ such
that $\dbar v = u$ and
$$
  \|v\|^2  \leq  \int_W \langle A^{-1}_{E,\omega} u, u \rangle dV.
$$
\end{izr}

\begin{op} If $v$ is replaced by the minimal $L^2$-norm solution and $u$ is smooth, so is $v.$
\end{op}
\medskip
The positivity of $A_{E,\omega}$ can be expressed with the coefficients of $i\Theta(E).$
In bidegree $(n,q)$  the positivity of the operator $A_{E,\omega}$ follows from Nakano positivity of $E.$ They are connected by the following formula with respect to the standard K\"{a}hler metric and an orthonormal frame on $E$ at a given point:
$$
  \langle A_{E,\omega} u, u \rangle = \sum_{|S| = q-1} \sum_{j,k,\lambda,\mu} c_{jk\lambda\mu}u_{jS,\lambda}\ol{u}_{kS,\mu}, \;\; u = \sum u_{J\lambda}dz_J \otimes e_{\lambda},
$$
where $i\Theta(E)$ is given by (\ref{tensor}).

In  bidegree  $(n,q)$ we have a theorem that provides the estimates in possibly noncomplete K\"{a}hler metric provided that
the manifold possesses a complete one. 

\begin{izr}[\rm Theorem VIII-6.1, \cite{dem}] \label{dem2} Let $(W, \hat{\omega})$ be a complete $n$-dimensional  K\"{a}hler manifold, $\omega$ another K\"{a}hler metric,
possibly non complete, and $E \ra W$ a Nakano semi-positive vector bundle. Let $u \in L^2_{n,q}(W,E),$ $q \geq 1,$
be a closed form satisfying
$$
  \int_W \langle A^{-1}_{E,\omega} u, u \rangle dV_{\omega} < \infty.
$$
Then  there exists $v \in L^2_{p,q-1}(W,E)$ such
that $\dbar v = u$ and

$$
  \|v\|^2\leq  \int_W \langle A^{-1}_{E,\omega} u, u\rangle dV_{\omega}.
$$
\end{izr}


\subsection{$\dbar$-equation in bidegree $(n,q)$}

\medskip

We can now solve the $\dbar$-problem for $(n,q)$-forms with the metric $h$ given by  theorem \ref{main thm1}. The  curvature tensor equals
$$
  i {\Theta}(E) = i\Theta(E)_0 + i \partial \dbar ((C+3)\Phi + \log \Phi_1)
$$
and therefore the curvature form $ A_{E,\omega}$  is strictly positive on the neighbourhood $\tilde{V}$ of $a(\ol{U} \setminus S),$ conic
along $a(S).$
Given $g : X \ra \Cc,$ $N(g) \supset S$ there exist by \cite{pre} an arbitrarily thin and sharp Stein neighbourhood $V \subset \tilde{V}$
of $a(U\setminus N(g))$ in $Z,$ conic along $a(N(g))$ and it possesses a complete K\"{a}hler
metric. As a result  theorem \ref{dem2} yields the following

\begin{figure}[h!]
\begin{center}
 \epsfysize=45mm
 \epsfbox{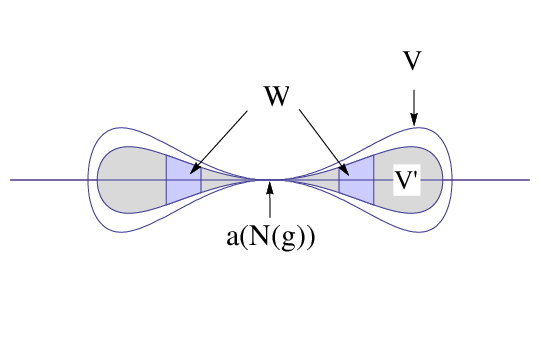}
 \caption{The cone $V'$ and the segment $W$}
 \label{segment}
\end{center}
  \end{figure}

\begin{izr}\label{rast} Let $u$ be a closed smooth $(n,q)$-form on $V$ with values in $E$ satisfying
$$
  \int_V \langle A^{-1}_{E,\omega} u,u\rangle_{h}e^{-M\log|g|} dV_{\omega} < \infty
$$
for some $M \geq 0.$ Then there exist a smooth $(n,q-1)$-form $v$ solving $\dbar v = u$ with
$$
  \|v\|^2 =   \int_V \langle v,v\rangle_{h}e^{-M\log|g|} dV_{\omega} \leq \int_V \langle A^{-1}_{E,\omega} u,u\rangle_{h}e^{-M\log|g|} dV_{\omega}.
$$
Assume in addition that  $q = 1$   and that the  smooth form $u$ has at most polynomial growth when approaching the boundary with respect to $h_Z$ and $h_E.$  Then $v$ has at most polynomial growth at the boundary. If $\|u\|_{\infty}$ is bounded and $M$ is large enough, then within a smaller cone
$V' \subset V$ with $\partial V \cap \partial V' \subset a(N(g))$ obtained by shrinking $V$ in the vertical direction and taking a smaller neighbourhood of $a(\ol{U})$ (see figure \ref{segment}) we have
$\lim_{z \ra z_0}v(z) = 0$ for every point $z_0 \in a(N(g)).$
\end{izr}

Notice that by multiplying the metric by $e^{-M \log |g|}$  we do not change the curvature, since $\log|g|$ is pluriharmonic.

The last statement of the theorem follows from Bochner-Martinelli-Koppelman (BMK) formula. Let $v$ be a $(p,0)$-form, $v(z) = \sum_{|P|=p} a_P(z) dz_P,$  and define
 $|v(z)|_{\infty}:= \max_P|a_P(z)|,$ $P$ is a multiindex. Rephrasing the proof in \cite{fl}, lemma 3.2., for $(p,0)$-forms we obtain

\begin{lm}\label{BMK} Let $v$ be a $(p,0)$-form with coefficients in $ {\mathcal C}^1(\ve B^n(0,1)),$ where $B^n(0,1)$ is the unit ball in $\Cc^n.$
Then we have  the estimate
$$
  |v(0)|_{\infty} \leq C(\ve^{-n}\|v\|_{L^2(\ve B^n(0,1))} + \ve \|\dbar v\|_{L^{\infty}(\ve B^n(0,1))}).
$$
The constant $C$ depends on $n$ only.
\end{lm}

\dok
Let $\chi$ be  a smooth cut-off function on $B = B^n(0,1),$ $\chi = 1$ on $\frac{1}{2} B.$
Fix a multiindex $P$ and estimate  $v(\zeta)_P = a(\zeta)_P d\zeta_P.$
The BMK kernel is
$$
  B(z,\zeta) =\frac{(n-1)!}{(2i\pi)^n |\zeta - z|^{2n} } \sum (-1)^{j-1} (\ol{\zeta}_j - \ol{z}_j) \wedge d(\ol{\zeta} - \ol{z})[j] \wedge d(\zeta -z),
$$
where $dz = dz_1 \wedge \ldots \wedge dz_n$ and $dz[j]$ is the $(n-1)$-form obtained from $dz$ by omitting $dz_j.$

We set $B = \sum B^p_q$ where $B^p_q$ is of the type $(p,q) $ in $z$ and $(n-p, n-q-1)$ in $\zeta$
and let $B^p_0 = \sum B^{p,P}_0$ where $B^{p,P}_{0}$  is of the type  $dz_P$.

The BMK formula gives
\begin{eqnarray*}
  (-1)^p v(0)_P &=& \int_{\partial \ve B} v(\zeta) \chi( \zeta/\ve) \wedge B^{p,P}_0(0,\zeta) -
                \int_{\ve B} \dbar(v(\zeta) \chi( \zeta/\ve)) \wedge B^{p,P}_0(0,\zeta)\\
              &=&- \int_{\ve B} \dbar v(\zeta) \wedge \chi(\zeta / \ve) B^{p,P}_0(0,\zeta) - \int_{\ve B} v(\zeta)\wedge \dbar(\chi(\zeta/\ve)) \wedge B^{p,P}_0(0,\zeta).\\
\end{eqnarray*}

In the second integral the form $\dbar\chi( \zeta/ \ve) \wedge B^{p,P}_0(0,\zeta)$ has its support on $\ve/2<|\zeta|< \ve$ and is ${\mathcal C}^{\infty}, $ $B^p_0$ has coefficients bounded by $\|\ve\|^{-2n+1}$, $\dbar(\chi( \zeta/\ve)) =
\dbar\chi(z)|_{z =  \zeta/\ve} \ve^{-1}$  and by Cauchy-Schwarz inequality the integral can be estimated by $\ve^{-n}C_1\|v\|_{L^2(\ve B)}.$ The first integral is bounded  by
$\ve C_2\|\dbar v\|_{L^{\infty}}.$\qed\\

\noi{\bf Proof of \ref{rast}.} The only thing to be proved is the last paragraph of the theorem. We have to compare the $L^2$-estimates for metric $h$ on $E$  and K\"{a}hler  form $\omega$ with
analogous estimates for the ambient Hermitian metric $h_E$ on $E$ and the ambient Hermitian form $\omega_Z$ on $Z.$ The distances on $Z$ are measured with respect to $h_Z.$ Near $a(S)$
we have $h(z,w) \approx \|z_2\|^{-\kappa_1}h_E(z,w),$ $\kappa_1 \in \Nn$  and
$dV_{\omega}(z,w) \approx \|z_2\|^{\kappa_2}dV_{\omega_Z}(z,w),$ $\kappa_2 \in \Nn;$  outside an open set containing $a(S)$ both metrics and both volume forms are uniformly equivalent. Therefore the result follows immediately for those boundary points which are not in $a(N(g)).$ Let  $\|v\|^2 = \|v\|^2_{(V,h |g(z)|^{-M},\omega)}.$
For small ball of radius $\delta/2$ and centre $z_0$ at the distance $\delta$ from $\partial V$ we have  the estimate
$$
   \|v\|^2  \geq \|v\|^2_{(B(z_0, \delta/2),h |g(z)|^{-M},\omega)} \geq  \inf_{B(z_0,\delta/2)}\frac{1}{|g(z)|^{M}} \|v\|^2_{(B(z_0, \delta/2),h,\omega)}.
$$
Let $z_0 = (z,w)$ be a point near $a(S).$ Then  we can estimate
\begin{eqnarray*}
   \|v\|^2_{(B(z_0, \delta/2),h_E,\omega_Z)} &\leq& \|v\|^2_{(B(z_0, \delta/2),h|g(z)|^{-M},\omega)} \sup_{B(z_0,\delta/2)}|g(z)|^{M}\|z_2\|^{\kappa} \\
    &\leq& \|v\|^2 \sup_{B(z_0,\delta/2)}{|g(z)|^{M}} \|z_2\|^{\kappa}
\end{eqnarray*}
for  $\kappa =\kappa_1 - \kappa_2 \in \Zz;$ the zeroes of the form $\omega$ may not be compensated by the poles of $h$ and so the
exponent $\kappa$ can be negative.

Near points in $a(N(g))$ we can estimate the sup norm of $v$ in the following way. Let $V' \subset V$ be a cone inside $V'$ as in theorem \ref{rast}
 (figure \ref{segment}).  The form $v$ is continuous on $\ol{V'} \setminus a(N(g)).$

Consider the segment $W = W(\ve):= \lbrace(z,w) \in \ol{V'}, \ve \leq |g(z)| < 2\ve \rbrace.$
The distance $\delta:= d(W(\ve), \partial V)$ with respect to $h_Z$
depends polynomially on $\ve$ and therefore polynomially on $|g|$ so together with the above lemma we conclude from the estimate
$$
  |v(z_0)|_{\infty} \leq C\left(\left(\frac{2}{\delta}\right)^{n}\|v\|  \sup_{B(z_0,\delta/2)}{|g(z)|^{M}}\|z_2\|^{\kappa} + \frac{\delta}{2} \|u\|_{\infty}\right)
$$
that the values of $v$ have at most polynomial poles when $z_0$ approaches $a(N(g)).$ If $\|u\|_{\infty}$ is bounded and   $M$ is large enough such that
$$\sup_{z_0\in V, d(z_0,W(\ve))< \delta/2}{|g(z)|^{M}}\|z_2\|^{\kappa} \ra 0$$ as $\ve \ra 0$
the values of $v$ go to $0$ when approaching ${a(N(g))}$ within $\ol{V'}.$ In this case $v$ has a continuous extension to $\ol{V'}.$

$\;$ \qed

Consider a neighbourhood $V_{z_0} \subset  Z$ of a point $z_0 \in a(S)$ with the standard K\"{a}hler metric $\omega_0 =
i\sum dz_i \wedge d\ol{z}_i.$ Let
$A_{E,\omega_0}$ be a commutator with respect to the standard metric $\omega_0$ in the bidegree $(n,q),$
$$
  A_{E,\omega_0} u = \sum_{|I| = q-1,j,k,\lambda,\mu} c_{jk\lambda\mu} u_{jI}dz_1\wedge \ldots  \wedge dz_n  \wedge d\ol{z}_{kI} \otimes e_{\mu}
$$
Then the largest eigenvalue of $A_{E,\omega_0}^{-1}$  has at most polynomial poles on $a(S).$
The commutator $A_{E,\omega}$ with respect to the given K\"{a}hler metric
$\omega = i\sum \gamma_i dz_i \wedge d\ol{z}_i$ is given by Lemma VIII-6.3,\cite{dem}:
$$
  A_{E,\omega} u = \sum_{|I| = q-1,j,k,\lambda,\mu} \gamma_j^{-1}c_{jk\lambda\mu} u_{jI}dz_1 \wedge \ldots \wedge dz_n \wedge d\ol{z}_{kI} \otimes e_{\mu}.
$$
Then the maximal eigenvalue of $A_{E,\omega}^{-1}$   still has at most polynomial poles on $Z_S,$ i.e. it behaves in the worst case as
$\|z_2\|^{-k}$ for some $k.$

\begin{pos}[Extensions]\label{posext} Notation as above. Let  $v$
be a section of $\Lambda^{n,0}T^*Z \otimes E|_{a(X)}$ defined on a neighbourhood of $a(\ol{U})$ , i.e.
 a  holomorphic $(n,0)$-form with values in $E$ and coefficients in ${\mathcal J}(S)^k.$    There exist  a Stein neighbourhood $V$ of $a(U \setminus N(g))$ conic along $a(N(g)) $ and a $(n,0)$-form $\tilde{v}  \in \Lambda^{n,0}T^*Z \otimes E|_{V}$ extending $v$ with at most polynomial growth at the boundary.
\end{pos}

\noi{\bf Proof of \ref{posext}}. Recall that $r_0 = \dim VT(Z)$  is the fibre dimension. Since $V$ is Stein it is  K\"{a}hler and complete and because $a(X)$ in $V$ is given as a zero set of finitely many global functions,  the K\"{a}hler manifold ${V \setminus a(X)}$ is also complete (lemma VIII-7.2, \cite{dem}; because $V$ is Stein the analytic set $a(X) \cap V$ is defined by finitely many holomorphic functions and then the bundle $E$ in lemma 7.2 is trivial). The function
$$
    \Phi_2 = \vph + \vph_1 + \log(\vph_1)
$$
is strictly plurisubharmonic on some conic neighbourhood of the form $\|w\| \leq \|z_2\|^{k_4}$ and has a logarithmic pole on the  section $a(U).$ This follows immediately from the estimates derived in the proof of  theorem \ref{main thm1}.
 We are solving the $\dbar$-equation with the metric
 $$
   h_4 = h e^{-r_0\Phi_2}.
 $$
   Take an extension of the form $v$ in the vertical direction obtained by patching together local holomorphic lifts as in remark \ref{vext}, denote it again by $v$  and let $u = \dbar v.$  Since $u(z,0) = 0$ close to $a(S)$ the coefficients of $u$ are bounded by $C \|w\| \|z_2\|^{k}$  and away from $a(S)$ by  $C\|w\|.$
 By construction we have $\vph_1 \geq \|w\|^2 \|z_2\|^{2 k_1}.$ The inverse of $A_{E,\omega}$ with respect to the metric $h_4$ has a polynomial pole on $a(S)$ and the metric $h_4$ has a polynomial pole there, so we have a polynomial pole in the scalar product. Let the whole term be bounded by $\|z_2\|^{-2 k_3}.$

Let us introduce the polar coordinates in the base and fibre directions in the integral
$$
  \int_{V \setminus a(X)} \langle A^{-1}_{E,\omega} u,u\rangle_{h} e^{-r_0 \Phi_2} dV_{\omega}.
$$
On a neighbourhood of  a point in  $a(S)$ the integrand is of the form
\begin{eqnarray*}
  &&(\|z_2\|^{-2 k_3}) (\|w\|^{2 }\|z_2\|^{2k}) (\|w\|^{-2 r_0 } \|z_2\|^{-2 r_0 k_1}) (\|w\|^{2r_0 - 1} \|z_2\|^{2k_5}\|z\|^{2\dim X - 1}) = \\
  && =  \|w\| \|z_2\|^{2(k - k_3 - r_0 k_1 + 2 k_5)}\|z\|^{2\dim X - 1}.
\end{eqnarray*}
The  terms in the last bracket come from the volume form if we introduce the polar coordinates in the base and fibre directions and take into account the form $\omega$ which has zeroes on $Z_S.$
The integral on some neighbourhood of this point is reduced to
\begin{eqnarray*}
  && c_1\int_0^{\delta}\|z\|^{2\dim X - 1}  d\|z\| \int_0^{\|z_2\|^{k_4}}  \|w\| \|z_2\|^{2(k - k_3 - r_0 k_1 + 2 k_5)} d\|w\|\leq \\
  && \leq c_2 \int_0^{\delta}\|z_2\|^{2(k - k_3 - r_0 k_1 + 2 k_5 +  k_4)}d\|z\|
\end{eqnarray*}
and it converges if either the cone is sharp enough (i.e. $k_4$ large) or the form has a zero of high enough order ($k$ large).

On a neighbourhood of points in $a(U \setminus S)$ the integral is approximately of the type $\|w\|$ and is therefore finite because the set $V$ is relatively compact.

Let $\tilde u$ be the solution of $\dbar \tilde u = u$ given by theorem \ref{rast}. The integrability condition
\begin{equation}\label{l2ocena}
  \|\tilde u\|^2_{V \setminus a(X)} =  \int_{V \setminus a(X)} \langle \tilde u ,\tilde u  \rangle_{h_3} e^{-r_0 \Phi_2} dV_{\omega} < \infty
\end{equation}
implies that on a neighbourhood of $z_0 \in a(U \setminus S),$ where $E$ is trivial, the forms  $dV_{\omega}$ and  $dV_{\omega_Z}$ are equivalent and $h_3$ is equivalent to $h_E,$  the section $\tilde u $ is in $L^2_{loc},$ because $\Phi_2$ has zeroes on $a(U).$
  Therefore the components $\tilde u _i$ are in $L^2_{loc}$ and because $u$ is smooth  the solution $\dbar \tilde u =u$ holds in the distribution sense on $V,$ so the section $\tilde v := v - \tilde u$ is holomorphic in the distribution sense and by ellipticity it is smooth (compare \cite{dem}, section VIII-7).
  Therefore $\tilde u$ is also smooth. Because $r_0 = \codim_Za(X),$ the weight $e^{-r_0 \Phi_2}$ is not locally integrable and since the integral (\ref{l2ocena}) exists the section $\tilde u$ must be zero on $a(U\setminus S).$
  Because the weight $e^{-r_0\Phi_2}$ has poles on $Z_S \cup a(X)$ and is bounded from below on $V$ we also have
  $$
  \int_{V} \langle \tilde u ,\tilde u  \rangle_{h_3} dV_{\omega} < \infty.
  $$
  The polynomial behaviour at the boundary now follows from the estimates in the proof of theorem \ref{rast} for $M = 0.$
\qed

\begin{op} If we  replaced $r_0$ by $\tilde{r}_0 > r_0 $ then $\tilde u$ would have to have zeroes of higher order on $a(U \setminus S)$ to insure the integrability of (\ref{l2ocena}).
Similar ideas work for jet interpolation at one point (not in $a(S)$), since we have a local holomorphic extension. The weight is defined as $M\log (\|z\|^2 + \|w\|^2)$ on a neighbourhood of the given point and continued as a constant outside. The negativity of the curvature created by such weight can be compensated by $e^{-c\Phi}$ since we are away from $a(S).$
\end{op}

\subsection{$\dbar$-equation in bidegree $(0,q)$}

In this section we prove a  theorem  analogous to theorem \ref{rast} for $(0,q)$-forms. In this case the positivity of the curvature tensor is no longer ensured by the positivity of the bundle curvature. Therefore  a $(0,q)$-form is viewed as a $(n,q)$-form with values in a different vector bundle.

Let the notation be as usual. Let $u \in  \Lambda^{0,q}T^*Z \otimes E|_{V'}$ where $V'\subset Z$ possesses a complete K\"{a}hler metric.
Let $h_{\omega}$ be the metric on $TZ$ induced by the  K\"{a}hler metric $\omega.$  The canonical pairing locally gives a decomposition
$1 = v\otimes v^*,$ where $v \in \Lambda^{n,0}T^*Z$ and $v^* \in \Lambda^{n,0}TZ.$
Thus $u$ can be viewed as a
$(n,q)$-form $\tilde u$ with values in $\tilde{E} = \Lambda^{n,0}TZ \otimes E.$ This adds an additional term to the curvature tensor, namely
the curvature of the determinant bundle $\det TZ = \Lambda^{n,0}TZ$ with respect to $h_{\omega}.$ The curvature is the Ricci curvature and so
the curvature tensor equals
$$
  i\Theta(\tilde{E}) = i \Id_{\det TZ }\otimes \Theta(E)  + \Ricci(\omega)\otimes \Id_E.
$$
 Assume that $E$ is trivial with local frame
$e_1,\ldots,e_{r}.$ In local coordinates $\zeta$ we have
\begin{eqnarray*}
  u = u^{\zeta} &=& \sum u^{\zeta}_{Q,\lambda} d\ol{\zeta}_{Q} \otimes e_{\lambda},  \\
\tilde{u} = \tilde{u}^{\zeta} &=& \sum u^{\zeta}_{Q,\lambda} d\ol{\zeta}_{Q} \wedge  d\zeta_1  \wedge \ldots \wedge d\zeta_n \otimes
  (\partial/\partial \zeta_1) \wedge\ldots \wedge (\partial/\partial \zeta_n) \otimes e_{\lambda}
\end{eqnarray*}
for multiindices $|Q| = q.$
Therefore $\tilde{u}$ is a form with values in $\tilde{E}.$ If $H_{\omega}$ is a  matrix representing $h_{\omega}$ and $h_{\omega}^*$ is the induced metric on the dual
\begin{eqnarray*}
  |\tilde{u}|^2(\zeta) &=& \sum u^{\zeta}_{Q,\lambda}(\zeta) \ol{u^{\zeta}}_{Q',\lambda'}(\zeta) \langle d\ol{\zeta}_{Q}, d\ol{\zeta}_{Q'} \rangle_{h_{\omega}}
       \cdot   \|d\zeta_1 \wedge \ldots \wedge  d\zeta_n \|^2_{h_{\omega}^*}\cdot \\
         && \cdot \|(\partial/\partial \zeta_1) \wedge \ldots \wedge  (\partial/\partial \zeta_n)\|^2_{h_{\omega}} \langle e_{\lambda}, e_{\lambda'}\rangle_{h}.
\end{eqnarray*}
Because $\|d\zeta_1 \wedge \ldots \wedge  d\zeta_n \|^2_{h_{\omega}^*} = \det(H_{\omega}^{-T})$ and  $\|(\partial/\partial \zeta_1) \wedge \ldots \wedge  (\partial/\partial \zeta_n)\|^2_{h_{\omega}} = \det H_{\omega}$ the norm is equal to the norm of $u.$

We would like to find a weight which removes the Ricci curvature. By proposition \ref{sect} with $E = \det TZ$ there exist finitely many
almost holomorphic sections $v_i,$ holomorphic  to order $l_3$ in $w$ with zeroes of order $k_3$ on $Z_S$  generating the $\det TZ$ away from $Z_S.$
The metric on the determinant bundle $h_{\det TZ}$ induced by $\omega$ defines the squares of the norms
$$
  f_i(z,w)= \langle v_i(z,w),v_i(z,w)\rangle_{h_{\det TZ}}.
$$
The function
$$
  \vph_2(z,w) = \sum \langle v_i(z,w),v_i(z,w) \rangle_{h_{\det TZ}}
$$
is defined on a neighbourhood $V_T$ of $a(\ol{U})$ and has locally  polynomial zeroes over $Z_S$ - the metric itself has polynomial zeroes and the vector fields have polynomial zeroes.

Let $v$ be a nonzero holomorphic section of the determinant bundle defined on a neighbourhood of a point $(z,0) \in a(S).$ Then the metric $h_{\det TZ}$  can be represented as multiplication by the function $f(z,w) =  \langle v(z,w),v(z,w) \rangle_{h_{\det TZ}}$ and the Ricci curvature equals $-i\partial\ol{\partial} \log f \Id_{\tilde{E}}.$

By construction we have $v_i = \va_i v$ for some functions $\va_i,$ holomorphic in the fibre direction, holomorphic to the degree $l_3$ with zeroes of order $k_3$ on the fibres over $a(S).$ This implies that
$$
  \vph_2 = \sum \langle v_i,v_i\rangle_{h_{\det TZ}} = \sum \va_i \ol{\va_i}\langle v,v \rangle_{h_{\det TZ}} = ( \sum |\va_i|^2) f = \|\va\|^2 f,
$$
where $\va$ is a vector with components $\va_i.$  The function $\|\va\|^2$ has zeroes only on $Z_S$ so we have the estimate $ \|\va\|^2 \geq c \|z_2\|^{2 k_3}.$
Let's multiply the metric $h$ by the weight
$$
  e^{-\log \vph_2}.
$$
The weight adds the term $(i\partial\ol{\partial} \log f + i\partial\ol{\partial} \log \|\va\|^2)\Id_{\tilde{E}}$ to the curvature thus killing the Ricci curvature
and adding a term which has bounded negative part in a conic neighbourhood (calculation is the same as in the section $3$).
As before we can  compensate the negativity of the curvature by multiplying the metric by the weight
$e^{-c\Phi}$ and at the same time achieve that the lowest eigenvalue decreases at most polynomially.
Denote the new metric by $h_5,$
$$
  h_5 = h e^{-(c\Phi + \log \vph_2)}.
$$
As a result for some large constant $c$ the curvature tensor with respect to $h_5$
$$
  i\Theta(\tilde{E}) = i\Id_{\det TZ}\otimes \Theta(E)  + \Ricci(\omega)\otimes \Id_E +
  i\partial\dbar(\log \vph_2 +  c \Phi)\otimes \Id_{\tilde{E}}
$$
is positive and this enables us to solve the $\dbar$-equation with at most polynomial growth at the boundary and with zeroes on $a(N(g)).$ If we view $(0,q)$-form $u$ as a $(n,q)$-form  we obtain as a corollary  to  theorem \ref{rast} the following
\begin{izr}\label{0qforme} Let $u$ be a closed smooth $(0,q)$-form on $V$ with
$$
  \int_V \langle A^{-1}_{\tilde{E},\omega} u,u\rangle_{h_5}e^{-M\log|g|} dV_{\omega}
$$
for some $M \geq 0.$ Then there exist a smooth $(0,q-1)$-form $v$ solving $\dbar v = u$ with
$$
  \|v\|^2 =   \int_V \langle v,v\rangle_{h_5}e^{-M\log|g|} dV_{\omega} \leq \int_V \langle A^{-1}_{\tilde{E},\omega} u,u\rangle_{h_5}e^{-M\log|g|} dV_{\omega}.
$$

\end{izr}

\begin{op}
Note that the sign of the Ricci curvature does not play any role since we are removing the Ricci curvature by the weight in contrast with the previous theorem where we needed the positivity of the induced curvature on the determinant bundle  in order to compensate for the possible negativity of the Hermitian metric.
\end{op}

\section{Vertical sprays on conic neighbourhoods}

This section is devoted to the proof of  theorem \ref{main thm2}. Consider the set $U.$  We are looking for sections which are defined on a conic neighbourhood of a given compact  set $a(\ol{U})$ and such that they generate the vertical tangent bundle  $VT(Z)$ on an open neighbourhood of $a(K).$
To avoid too many notations we use the letter $U$ for such a neighbourhood  and will  shrink $U$ if necessary. Let ${\mathcal{VT}}(Z)$ denote the sheaf of sections of  $VT(Z).$ Let $v_i$ be
almost holomorphic sections of $VT(Z),$ holomorphic to the degree $l_4$ in $w$ and with zeroes  of order $k_4$ given by
proposition \ref{sect}. Let $u_i =\dbar v_i$ and view it as a $(n,1)$-form as in the previous section.
 Define the metric
 $$
    h_6 = h_5 e^{-r_1\Phi_2}.
 $$
 We have to show that over a suitable conic neighbourhood $V_1$ the integral
$$
  I = \int_{V_1\setminus a(X)} \langle A^{-1}_{\tilde{E},\omega} u_i,u_i\rangle_{h_5} e^{-r_1 \Phi_2} dV_{\omega}
$$
is convergent for $r_1 \geq r_0;$ recall that $r_0$ is the fibre dimension. The integrability is problematic only on  neighbourhoods of  points in $a(U).$ Let us first consider points in $a(S).$ The terms in the integrand are of the following form:  the form $u_i$
is of the type $\|w\|^{l_4+1}\|z_2\|^{k_4}$  and $A^{-1}_{\tilde{E},\omega}$ and ${h_5}$ have in the worst case  polynomial poles in $\|z_2\|.$ Let the scalar product
 $\langle A^{-1}_{\tilde{E},\omega} u_i,u_i\rangle_{h_5}$ be of the form $\|w\|^{2l_4+2}\|z_2\|^{2k_4 - n_1}.$
 The weight $e^{-r_1\log\Phi_2 - \log \vph_2}$ has the type $\|w\|^{-2r_1}\|z_2\|^{-2k_1r_1 - 2k_3}$ and
$dV_{\omega}$ is of the type $\|z_2\|^{2k_5}dV_{h_Z}.$  After introducing the polar coordinates in $w$ and $z$ on a neighbourhood $V_{z_0}$ of the point   $z_0 \in a(S)$ the integral
$I_1 = \int_{V_0 \cap (V_1\setminus a(X))} \langle A^{-1}_{\tilde{E},\omega} u_i,u_i\rangle_{h_5} e^{-r_1 \Phi_2} dV_{\omega}$
takes the form
\begin{eqnarray*}
  I_1 \leq &&\mbox{\rm const} \int_0^{\delta}\|z_2\|^{-n_1 + 2 k_4 - 2 r_1 k_1 -2k_3+ 2k_5}\|z\|^{2 \dim X - 1}   \cdot\\
  &&\cdot \int_0^{\|z_2\|^{k_6}}\|w\|^{2 l_4 + 2 - 2 (r_1 - r_0) - 1} d\|w\|\,  d\|z\|,
\end{eqnarray*}
where $\|w\| \leq \|z_2\|^{k_6}$ describes the type of the cone near $a(S).$

Put $r_1 = r_0.$ Then, if either  $k_4$ is large, meaning that the initial vector fields have zeroes of high order on $a(S),$ or the cone is sharp enough, for example $k_6 > n_1,$ or the vector fields are holomorphic to a very high order ($l_4$ large) the integral converges.
Near points from $a(U \setminus S)$ we only have the inner integral with $\|z_2\|^{k_6}$ replaced by some fixed $\delta$  and it converges for $l_4 \geq 0.$
Even if we start with any vector field  with zeroes of high order on $a(S)$ and construct an extension $v$ by remark \ref{vext} the integral converges. In this case we have $l_4 = 0.$

If $r_1 > r_0,$ then again  near $a(U \setminus S)$ the  integral converges if $l_4 > r_1 - r_0.$
Note that $r_1 - r_0$ is approximately the order of the jet interpolation and if the result is supposed to give a holomorphic section then the initial section must already be holomorphic to a high degree.

Thus there exist a neighbourhood $V_1$ of $a(U \setminus S)$, conic along $a(N(g)$ and such that $I < \infty.$ Then  theorem \ref{0qforme} for $q = 1, M=0$ with $h_6$ instead of $h_5$ and $V_1 \setminus a(X)$ instead of $V$  yields
the vector fields $\tilde{v}_i$ with values in $ VT(Z)$ of polynomial growth at the boundary. As in the proof of corollary \ref{posext} we show that the sections $\tilde{v}_i$ are zero on $a(U\setminus S).$ Moreover,
we have
$$
  \int_{V_1} \langle \tilde{v},\tilde{v}\rangle_{h_5} dV_{\omega} < \infty.
$$
The holomorphic vector fields  $v_i - \tilde{v_i}$  still generate the $VT(Z)$ on a  neighbourhood of $a(U \setminus S).$
In particular, they generate the bundle on a neighbourhood of $a(K)$ in $Z.$

We have to show that the vector fields can be corrected to vector fields with zeroes on $a(N(g)).$
If we take a slightly thinner and sharper cone along $a(N(g))$
and shrink $U$ a little they will be bounded when away from $a(N(g))$.  Denote this conic neighbourhood by $V'_1.$ As in the proof of \ref{rast} we see that if we approach $a(N(g)$  within $V'_1$   the vector fields $\tilde{v}_i$  have at most polynomial poles on $a(N(g)).$

But then the vector fields
$g^k(v_i -\tilde{v_i}) $ for sufficiently large $k$ still generate the bundle wherever $v_i - \tilde{v}_i$ do and approach $0$  as $|g| \ra 0$ as fast as we want. In particular they are (at least) continuous on the closure of $V'_1$ and can be extended to global continuous vector fields with zeroes on $Z_{N(g)}.$
Let $V$ be a smaller conic Stein neighbourhood inside $V'_1.$
The flows  $\vph_{i,t_i}(z)$ of the fields of $v_i$  remain in $V$ for $z$ in a thinner and sharper conic neighbourhood $V'$ (see figure \ref{okolice}) for sufficiently small times   and so generate a continuous vertical spray $s:=\vph_{1,\cdot} \circ \ldots \circ \vph_{m,\cdot} : V_Z \times B^m(0,\ve) \ra Z,$ $s(z,(t_1,\ldots,t_m)):= \vph_{1,t_1} \circ \ldots \circ \vph_{m,t_m}(z)$ on
 sufficiently small neighbourhood $V_Z$ of $a(U)$ in $Z,$ such that $s(z,t) \in V$ for $z \in V'.$ The restriction of $s$ to $a(U) \times B^m(0,\ve)$ is smooth and holomorphic on $a(U\setminus N(g)) \times B^m(0,\ve)$ and is therefore holomorphic on $a(U) \times B^m(0,\ve)$ since $N(g)$ is analytic.
This completes the proof of the main theorem in \cite{pre} in the case of manifolds.

\section*{Acknowledgments}
{The author was supported by  Slovenian Research Agency,  Research program P1-0291 and Research project J1-5432. The part of this paper was written while the author was visiting the NTNU, Trondheim, Norway and she wishes to thank this institution for its hospitality. The author also wishes to thank the reviewer for the comments which helped to improve the exposition.}

\end{document}